\documentclass[11pt,reqno,a4paper]{amsart}
\usepackage{amssymb}
\input amssym.def
\usepackage{amsmath,amsfonts,xcolor,mathtools,booktabs}
\usepackage[bookmarks=false]{hyperref}
\usepackage{amscd}
\usepackage[mathscr]{eucal}
\usepackage{enumitem}
\usepackage{orcidlink}
\allowdisplaybreaks
\numberwithin{equation}{section}

\hypersetup{colorlinks=true,citecolor={purple},linkcolor={teal},urlcolor={violet}}

\theoremstyle{plain}
\newtheorem{theorem}{Theorem}[section]

\newtheorem{lemma}[theorem]{Lemma}

\theoremstyle{definition}

\newtheorem{remark}[theorem]{Remark}

\makeatletter
\@namedef{subjclassname@2020}{%
\textup{2020} Mathematics Subject Classification}
\makeatother

\subjclass[2020]{11F03, 11F11, 33D15}

\keywords{Gosper-type identities, Lambert series, $\eta$-quotients, generalized $\eta$-quotients}

\begin{document}
\title[New Gosper-type Lambert series identities of levels $12$ and $16$]{New Gosper-type Lambert series identities of levels $12$ and $16$}	
\author[Russelle Guadalupe]{Russelle Guadalupe\orcidlink{0009-0001-8974-4502}}
\address{Institute of Mathematics, University of the Philippines Diliman\\
Quezon City 1101, Philippines}
\email{rguadalupe@math.upd.edu.ph}

\begin{abstract}
We derive new Gosper-type Lambert series identities of levels $12$ and $16$ using certain sums of generalized $\eta$-quotients on the genus zero congruence subgroups $\Gamma_0(12)$ and $\Gamma_0(16)$.
\end{abstract}

\maketitle

\section{Introduction}\label{sec1}
For an element $\tau$ in the complex upper half-plane $\mathbb{H}$, we define $q:=e^{2\pi i \tau}$, and for any complex number $a$, we define $(a;q)_{\infty} := \prod_{n=0}^\infty (1-aq^n)$ and 
\begin{align*}
(a_1,a_2,\ldots,a_m;q)_\infty:= (a_1;q)_\infty (a_2;q)_\infty\cdots (a_m;q)_\infty. 
\end{align*}
For any complex numbers $a$ and $b$ with $|ab| < 1$, the Ramanujan's general theta function is defined by 
\begin{align*}
f(a,b) =\sum_{n=-\infty}^\infty a^{n(n+1)/2}b^{n(n-1)/2} =(-a,-b,ab;ab)_\infty,
\end{align*}
where the last equality follows from the Jacobi triple product identity \cite[p. 35, Entry 19]{berndt3}. Recall that for a given sequence $\{a_n\}_{n\geq 1}$ of complex numbers, the Lambert series, introduced by Lambert \cite{lamb} in his work on the convergence of power series, is an infinite series of the form 
\begin{align*}
\sum_{n\geq1}\dfrac{a_nq^n}{1-q^n}.
\end{align*}
An example of a Lambert series identity is given by
\begin{align*}
\sum_{n\geq 1}\sigma(n)q^n=\sum_{n\geq 1}\dfrac{nq^n}{1-q^n} = \sum_{n\geq 1}\dfrac{q^n}{(1-q^n)^2},
\end{align*}
where $\sigma(n)$ is the sum of all positive divisors of $n$. In his study on the properties of the $q$-analogues of trigonometric functions, Gosper \cite{gosper} introduced the function 
\begin{align*}
\Pi_q := q^{1/4}\dfrac{(q^2;q^2)_\infty^2}{(q;q^2)_\infty^2}=q^{1/4}\psi^2(q),
\end{align*} 
where $\psi(q) := f(q,q^3)=(q^2;q^2)_\infty/(q;q^2)_\infty$, and conjectured the following Lambert series identities involving $\Pi_q$:
\begin{align}
&\sum_{n\geq 1}\dfrac{q^n}{(1-q^n)^2}-2\sum_{n\geq 1}\dfrac{q^{2n}}{(1-q^{2n})^2} = \dfrac{1}{24}\left(\dfrac{\Pi_q^4}{\Pi_{q^2}^2}-1\right)+\dfrac{2}{3}\Pi_{q^2}^2,\label{eq11}\\
&\dfrac{1}{\Pi_{q^5}^2}\left(\sum_{n\geq 1}\dfrac{q^{2n-1}}{(1-q^{2n-1})^2}-5\sum_{n\geq 1}\dfrac{q^{10n-5}}{(1-q^{10n-5})^2}\right)= \sqrt{\dfrac{\Pi_q^3}{\Pi_{q^5}^3}-2\dfrac{\Pi_q^2}{\Pi_{q^5}^2}+5 \dfrac{\Pi_q}{\Pi_{q^5}}},\label{eq12}\\
&6\left(\sum_{n\geq 1}\dfrac{q^n}{(1-q^n)^2}-5\sum_{n\geq 1}\dfrac{q^{5n}}{(1-q^{5n})^2}\right)+1\nonumber\\
&=\left(\dfrac{\Pi_q}{\Pi_{q^5}}+2+5\dfrac{\Pi_{q^5}}{\Pi_q}\right)\left(\sum_{n\geq 1}\dfrac{q^{2n-1}}{(1-q^{2n-1})^2}-5\sum_{n\geq 1}\dfrac{q^{10n-5}}{(1-q^{10n-5})^2}\right),\label{eq13}\\
&\dfrac{1}{\Pi_{q^6}^2}\left(\sum_{n\geq 1}\dfrac{q^{2n-1}}{(1-q^{2n-1})^2}-6\sum_{n\geq 1}\dfrac{q^{12n-6}}{(1-q^{12n-6})^2}\right)=\dfrac{\Pi_{q^2}^2}{\Pi_{q^6}^2}+2\dfrac{\Pi_{q^2}}{\Pi_{q^6}},\label{eq14}\\
&\dfrac{1}{\Pi_{q^9}^2}\left(\sum_{n\geq 1}\dfrac{q^{2n-1}}{(1-q^{2n-1})^2}-9\sum_{n\geq 1}\dfrac{q^{18n-9}}{(1-q^{18n-9})^2}\right)\nonumber\\
&=\left(\dfrac{\Pi_q}{\Pi_{q^9}}+3\right)\sqrt{\left(\dfrac{\Pi_q}{\Pi_{q^9}}\right)^{3/2}-3\dfrac{\Pi_q}{\Pi_{q^9}}+3\left(\dfrac{\Pi_q}{\Pi_{q^9}}\right)^{1/2}},\label{eq15}\\
&3\left(\sum_{n\geq 1}\dfrac{q^n}{(1-q^n)^2}-9\sum_{n\geq 1}\dfrac{q^{9n}}{(1-q^{9n})^2}\right)+1\nonumber\\
&=\left(\sqrt{\dfrac{\Pi_q}{\Pi_{q^9}}}+3\sqrt{\dfrac{\Pi_{q^9}}{\Pi_q}}\right)\left(\sum_{n\geq 1}\dfrac{q^{2n-1}}{(1-q^{2n-1})^2}-9\sum_{n\geq 1}\dfrac{q^{18n-9}}{(1-q^{18n-9})^2}\right).\label{eq16}
\end{align}
El Bachraoui \cite{bach} used the methods of Gosper and integrals of Ramanujan's theta functions to prove (\ref{eq11}). He \cite{he1,he2} confirmed (\ref{eq12}) and (\ref{eq13}) using Ramanujan's modular equations of degree five and proved (\ref{eq15}) via modular forms. Wang \cite{wang} observed that the above identities can be seen as equalities of certain modular forms on the congruence subgroup $\Gamma_0(N)$ and called such equality an identity of level $N$. Wang \cite{wang} also proved these identities and devised a method of discovering Lambert series identities of level $N$ via Eisenstein series, a basis for the vector space of modular forms on $\Gamma_0(N)$ of some positive weight, and the Sturm bound. Yathirajsharma, Harshitha, and Vasuki \cite{yhv} deduced the above identities except (\ref{eq14}) by applying Ramanujan's theta function identities. Recently, Yathirajsharma \cite{yat} established another form of (\ref{eq15}) by employing some $q$-trigonometric identities.

In this paper, we offer another proof of the identity (\ref{eq14}) of level $12$ and derive new Gosper-type Lambert series identities of levels $12$ and $16$ using sums of generalized $\eta$-quotients on $\Gamma_0(N)$ for $N\in\{12,16\}$, as shown in our main results below.

\begin{theorem}[level 12 identity]\label{thm11}
We have 
\begin{align}\label{eq17}
	&\dfrac{1}{\Pi_{q^6}^2}\left[24\left(\sum_{n\geq 1}\dfrac{q^{n}}{(1-q^{n})^2}-6\sum_{n\geq 1}\dfrac{q^{6n}}{(1-q^{6n})^2}\right)+5\right]\nonumber\\
	&\phantom{identities level shown below using}=5\dfrac{\Pi_{q^2}^3}{\Pi_{q^6}^3}+24\dfrac{\Pi_{q^2}^2}{\Pi_{q^6}^2}+42\dfrac{\Pi_{q^2}}{\Pi_{q^6}}+9\dfrac{\Pi_{q^6}}{\Pi_{q^2}}.
\end{align}
\end{theorem}

\begin{theorem}[level 16 identities]\label{thm12} 
We have 
\begin{align}
	&\dfrac{1}{\Pi_{q^8}^2}\left(\sum_{n\geq 1}\dfrac{q^{2n-1}}{(1-q^{2n-1})^2}-8\sum_{n\geq 1}\dfrac{q^{16n-8}}{(1-q^{16n-8})^2}\right)= 
	\dfrac{\Pi_{q^4}^3}{\Pi_{q^8}^3}+2\dfrac{\Pi_{q^4}^2}{\Pi_{q^8}^2}+4\dfrac{\Pi_{q^4}}{\Pi_{q^8}}+4,\label{eq18}\\
	&\dfrac{1}{\Pi_{q^8}^2}\left[24\left(\sum_{n\geq 1}\dfrac{q^{n}}{(1-q^{n})^2}-8\sum_{n\geq 1}\dfrac{q^{8n}}{(1-q^{8n})^2}\right)+7\right]\nonumber\\
	&\phantom{identities level shown below using}=7\dfrac{\Pi_{q^4}^4}{\Pi_{q^8}^4}+24\dfrac{\Pi_{q^4}^3}{\Pi_{q^8}^3}+72\dfrac{\Pi_{q^4}^2}{\Pi_{q^8}^2}+96\dfrac{\Pi_{q^4}}{\Pi_{q^8}}+112.\label{eq19}
\end{align}
\end{theorem}
Our main tool for proving (\ref{eq14}) and Theorems \ref{thm11} and \ref{thm12} is the following identity of Bailey \cite{bailey2} given by 
\begin{align}
\sum_{n=-\infty}^\infty\left[\dfrac{aq^n}{(1-aq^n)^2}-\dfrac{bq^n}{(1-bq^n)^2}\right]=a(q;q)_{\infty}^6\dfrac{f(-ab,-\frac{q}{ab})f(-\frac{b}{a},-\frac{aq}{b})}{f^2(-a,-\frac{q}{a})f^2(-b,-\frac{q}{b})}\label{eq110}
\end{align}
for $|ab| < 1$, which is a special case of his $_6\psi_6$ summation formula \cite[(4.7)]{bailey1}. This formula transforms a given Lambert series identity of level $N$ into sums of products of theta functions, which in turn can be written as sums of generalized $\eta$-quotients on $\Gamma_0(N)$. We then exploit the fact that $\Gamma_0(N)$ has genus zero for $N\in\{12,16\}$, which will be essential in establishing these identities by expressing them in terms of quotients of $\Pi_q$. We remark that our method presented in this paper is different from that of He \cite{he2} and Wang \cite{wang} in that we only apply modular functions instead of modular forms of positive weight and the Sturm bound to derive the above identities. In particular, the author \cite{guad} recently obtained two new Gosper-type Lambert series identities of level $14$ using this method.

We organize the rest of the paper as follows. We give a short background on modular functions, particularly $\eta$-quotients and generalized $\eta$-quotients, on the congruence subgroups $\Gamma_1(N)$ and $\Gamma_0(N)$ in Section \ref{sec2}. Applying Bailey's formula (\ref{eq110}) and some properties of sums of generalized $\eta$-quotients on $\Gamma_0(N)$ for $N\in\{12,16\}$, we establish identity (\ref{eq14}) and Theorem \ref{thm11} in Section \ref{sec3}, and Theorem \ref{thm12} in Section \ref{sec4}. In the proofs of these identities, we express the corresponding sum of generalized $\eta$-quotients on $\Gamma_0(N)$ in terms of the generator of the field of modular functions on $\Gamma_0(N)$ involving quotients of $\Pi_q$. We have executed most of our computations using \textit{Mathematica}. 

\section{Modular functions on $\Gamma_1(N)$ and $\Gamma_0(N)$}\label{sec2}

We consider in this paper modular functions on congruence subgroups
\begin{align*}
\Gamma_1(N) &:= \left\lbrace\begin{bmatrix}
	a & b\\ c& d
\end{bmatrix} \in \mathrm{SL}_2(\mathbb{Z}) : a-1\equiv c\equiv d-1\equiv 0\pmod N\right\rbrace,\\
\Gamma_0(N) &:= \left\lbrace\begin{bmatrix}
	a & b\\ c& d
\end{bmatrix} \in \mathrm{SL}_2(\mathbb{Z}) : c\equiv 0\pmod N\right\rbrace.
\end{align*}

We describe the action of an element of a congruence subgroup $\Gamma$ on the extended upper half-plane $\mathbb{H}^\ast := \mathbb{H}\cup \mathbb{Q}\cup \{\infty\}$ by a linear fractional transformation. We define the cusps of $\Gamma$ as the equivalence classes of $\mathbb{Q}\cup \{\infty\}$ under this action. We define a modular function on $\Gamma$ as a meromorphic function $f:\mathbb{H}^\ast\rightarrow\mathbb{C}$ such that $f(\gamma\tau)=f(\tau)$ for all $\gamma\in \Gamma$ and the $q$-expansion of $f$ at each cusp $r$ of $\Gamma$ is 
\begin{align*}
f(\gamma\tau) = \sum_{n\geq n_0} a_nq^{n/h}
\end{align*}
for some $\gamma\in \mathrm{SL}_2(\mathbb{Z})$ with $\gamma(\infty)=r$ and integers $h$ and $n_0$ with $a_{n_0}\neq 0$. We call $n_0$ the order of $f(\tau)$ at $r$, denoted by $\mathrm{ord}_r f(\tau)$, and $h$ the width of $r$, which is the smallest positive integer such that $\gamma[\begin{smallmatrix}
1 & h\\0 & 1
\end{smallmatrix}]\gamma^{-1}\in \pm\Gamma$. We say that $f(\tau)$ has a zero (resp., a pole) at $r$ if $\mathrm{ord}_r\,f(\tau) > 0$ (resp.,  $\mathrm{ord}_r\,f(\tau) < 0$).

Cho, Koo, and Park \cite{chokoop} characterized the inequivalent cusps of $\Gamma_0(N)$ and $\Gamma_1(N)$, and obtained their respective widths.

\begin{lemma}\label{lem21}
Let $a,a',c$ and $c'$ be integers with $\gcd(a,c)=\gcd(a',c')=1$. Let $S_{\Gamma_0(N)}$ be the set of all inequivalent cusps of $\Gamma_0(N)$. We denote $\pm1/0$ as $\infty$. Then 
\begin{enumerate}
	\item[(1)] $a/c$ and $a'/c'$ are equivalent over $\Gamma_0(N)$ if and only if there exist an integer $n$ and an element $\overline{s}\in (\mathbb{Z}/N\mathbb{Z})^\times$ such that $(a',c')\equiv(\overline{s}^{-1}a+nc,\overline{s}c)\pmod N$,
	\item[(2)] we have 
	\begin{align*}
		S_{\Gamma_0(N)} = \{a_{c,j}/c\in \mathbb{Q} : 0 < c \mid N, 0 < a_{c,j}\leq N, \gcd(a_{c,j},N)=1,\\ a_{c,j}= a_{c,j'}\stackrel{\text{def}}{\iff} a_{c,j}\equiv a_{c,j'}\pmod{\gcd(c,N/c)}\},
	\end{align*}
	\item[(3)] the width of the cusp $a/c\in S_{\Gamma_0(N)}$ is $N/\gcd(c^2,N)$.
\end{enumerate}
\end{lemma}

\begin{proof}
See \cite[Corollary 4(1)]{chokoop}.
\end{proof}

\begin{lemma}\label{lem22}
Let $a,a',c$ and $c'$ be integers with $\gcd(a,c)=\gcd(a',c')=1$. Let $S_{\Gamma_1(N)}$ be the set of all inequivalent cusps of $\Gamma_1(N)$. We denote $\pm1/0$ as $\infty$. Then 
\begin{enumerate}
	\item[(1)] $a/c$ and $a'/c'$ are equivalent under $\Gamma_1(N)$ if and only if there exists $n\in\mathbb{Z}$ such that $[\begin{smallmatrix}
		a' \\ c'
	\end{smallmatrix}]\equiv \pm[\begin{smallmatrix}
		a+nc \\ c
	\end{smallmatrix}]\pmod{N}$,
	\item[(2)] we have
	\begin{align*}
		S_{\Gamma_1(N)} =&\{y_{c,j}/x_{c,i}\in\mathbb{Q} : 0 < c\mid N, 0<s_{c,i},a_{c,j}\leq N,\\ 
		&\gcd(s_{c,i},N)=\gcd(a_{c,j},N)=1,\\
		& s_{c,i}=s_{c,i'}\stackrel{\text{def}}{\iff} s_{c,i}\equiv \pm s_{c,i'}\pmod{N/c},\\
		& a_{c,j}=a_{c,j'}\stackrel{\text{def}}{\iff} a_{c,j}\equiv \begin{cases}
			\pm a_{c,j'}\pmod{c} &\text{ if } c\in\{N/2,N\},\\
			a_{c,j'}\pmod{c} &\text{ otherwise, }
		\end{cases}\\
		&(x_{c,i},y_{c,j})\in\mathbb{Z}^2, \gcd(x_{c,i},y_{c,j})=1,\\
		&x_{c,i}\equiv cs_{c,i}\pmod{N}, y_{c,j}\equiv a_{c,j}\pmod{N}\},
	\end{align*}
	and the width of the cusp $a/c$ in $\Gamma_1(N)$ is $1$ (respectively, $N/\gcd(c,N)$) if $N=4$ and $\gcd(c,4)=2$ (respectively, otherwise).
\end{enumerate}
\end{lemma}

\begin{proof}
See \cite[Corollary 4(2)]{chokoop}.
\end{proof}

Let $\mathcal{M}^\infty(N)$ be the set of all modular functions on $\Gamma_0(N)$ that have poles only at infinity. Since $\Gamma_0(N)$ has genus zero for $N\in \{12,16\}$, we see that $\mathcal{M}^\infty(N)$ is generated by a single function $g\in \mathcal{M}^\infty(N)$. The following result shows that for these values of $N$, every element of $\mathcal{M}^\infty(N)$ can be expressed in terms of $g$. 

\begin{lemma}\label{lem23}
Let $N\in \{12,16\}$ and $h\in\mathcal{M}^\infty(N)$. Set $m := \mathrm{ord}_\infty\,h$. Then $h$ can be written as a polynomial in $g$ of degree $-m$ with complex coefficients. 
\end{lemma}

\begin{proof}
This follows from \cite[Theorem 8.3]{pauleradu}.
\end{proof}

We next define an $\eta$-quotient given by 
\begin{align*}
f(\tau) = \prod_{\delta\mid N} \eta^{r_{\delta}}(\delta\tau)
\end{align*}
for some indexed set $\{r_\delta\in\mathbb{Z} : \delta\mid N\}$, where $\eta(\tau)=q^{1/24}(q;q)_\infty$ is the Dedekind eta function. The following results give conditions for an $\eta$-quotient to be modular on $\Gamma_0(N)$. 

\begin{lemma}\label{lem24}
Let $f(\tau) = \prod_{\delta\mid N} \eta^{r_{\delta}}(\delta\tau)$ be an $\eta$-quotient with $\sum_{\delta\mid N} r_{\delta}=0$. Then $f$ is a modular function on $\Gamma_0(N)$ with character
\begin{align*}
	\chi(d):=\left(\dfrac{\prod_{\delta\mid N} \delta^{r_\delta}}{d}\right),
\end{align*}
where $(\frac{\cdot}{d})$ is the Kronecker symbol, if and only if
\begin{equation*}
	\sum_{\delta\mid N} \delta r_{\delta} \equiv 0\pmod{24}\quad\text{ and }\quad\sum_{\delta\mid N} \dfrac{N}{\delta}r_{\delta} \equiv 0\pmod {24}.
\end{equation*}
\end{lemma} 

\begin{proof}
See \cite[Proposition 5.9.2]{cohen}.
\end{proof}

\begin{lemma}\label{lem25}
Let $c, d$, and $N$ be positive integers with $d\mid N$ and $\gcd(c,d)=1$ and let $f(\tau) = \prod_{\delta\mid N} \eta^{r_{\delta}}(\delta\tau)$ be an $\eta$-quotient satisfying the conditions of Lemma \ref{lem24}. Then the order of $f(\tau)$ at the cusp $c/d$ is 
\begin{align*}
	\mathrm{ord}_{a/c}\,f(\tau) = \dfrac{N}{24d\gcd(d,\frac{N}{d})}\sum_{\delta\mid N} \gcd(d,\delta)^2\cdot\dfrac{r_{\delta}}{\delta}.
\end{align*}
\end{lemma}

\begin{proof}
See \cite[Lemma 4.1.3]{aygin}.
\end{proof}

We now define a generalized $\eta$-quotient given by
\begin{align*}
f(\tau) = \prod_{1\leq g\leq \lfloor N/2\rfloor}\eta_{N,g}^{r_g}(\tau)
\end{align*}
for some indexed set $\{r_g \in\mathbb{Z}: 1\leq g\leq \lfloor N/2\rfloor\}$, where
\begin{align*}
\eta_{N, g}(\tau)=q^{NB_2(g/N)/2}(q^g,q^{N-g};q^N)_\infty
\end{align*} 
is the generalized Dedekind eta function with $B_2(t) := t^2-t+1/6$. The following results provide the transformation formula of $\eta_{N, g}$ on $\Gamma_0(N)$ and conditions for a generalized $\eta$-quotient to be modular on $\Gamma_1(N)$.

\begin{lemma}\label{lem26}
The function $\eta_{N,g}(\tau)$ satisfies $\eta_{N,g+N}(\tau)=\eta_{N,-g}(\tau)=-\eta_{N,g}(\tau)$. Moreover, let $\gamma = [\begin{smallmatrix}
	a & b\\ cN & d
\end{smallmatrix}]\in \Gamma_0(N)$ with $c\neq 0$. Then we have 
\begin{align*}
	\eta_{N,g}(\gamma\tau) = \varepsilon(a,bN,c,d)e^{\pi i(g^2ab/N - gb)}\eta_{N,ag}(\tau),
\end{align*}
where 
\begin{align*}
	\varepsilon(a,b,c,d) = \begin{cases}
		e^{\pi i(bd(1-c^2)+c(a+d-3))/6}, & \text{ if }c\equiv 1\pmod{2},\\
		-ie^{\pi i(ac(1-d^2)+d(b-c+3))/6}, & \text{ if }c\equiv 0\pmod{2}.
	\end{cases}
\end{align*}
\end{lemma}

\begin{proof}
See \cite[Corollary 2]{yang1}.
\end{proof}

\begin{lemma}\label{lem27}
Suppose $f(\tau)=\prod_{1\leq g\leq \lfloor N/2\rfloor}\eta_{N,g}^{r_g}(\tau)$ is a generalized $\eta$-quotient such that
\begin{enumerate}
	\item[(1)] $\sum_{1\leq g\leq \lfloor N/2\rfloor}r_g\equiv 0\pmod{12}$, 
	\item[(2)] $\sum_{1\leq g\leq \lfloor N/2\rfloor}gr_g\equiv 0\pmod{2}$, and 
	\item[(3)] $\sum_{1\leq g\leq \lfloor N/2\rfloor}g^2r_g\equiv 0\pmod{2N}$.
\end{enumerate}
Then $f(\tau)$ is a modular function on $\Gamma_1(N)$. 
\end{lemma}

\begin{proof}
See \cite[Corollary 3]{yang1}.
\end{proof}

\begin{lemma}\label{lem28}
Let $N$ be a positive integer and $\gamma=[\begin{smallmatrix}
	a & b\\ c& d
\end{smallmatrix}]\in\mathrm{SL}_2(\mathbb{Z})$. Then the first term of the $q$-expansion of $\eta_{N, g}(\gamma\tau)$ is $\varepsilon q^{\delta}$, where $|\varepsilon|=1$ and 
\begin{align*}
	\delta = \dfrac{\gcd(c,N)^2}{2N}P_2\left(\dfrac{ag}{\gcd(c,N)}\right)
\end{align*}
with $P_2(t):=B_2(\{t\})$ and $\{t\}$ is the fractional part of $t$.
\end{lemma}

\begin{proof}
See \cite[Lemma 2]{yang1}.
\end{proof}

\begin{remark}\label{rem29}
We infer from Lemmas \ref{lem22} and \ref{lem28} that if $f(\tau)=\prod_{1\leq g\leq \lfloor N/2\rfloor}\eta_{N,g}^{r_g}(\tau)$ is a generalized $\eta$-quotient on $\Gamma_1(N)$ with $N > 4$, then the order of $f(\tau)$ at the cusp $a/c$ of $\Gamma_1(N)$ is 
\begin{align*}
	\mathrm{ord}_{a/c}\, f(\tau) = \dfrac{\gcd(c,N)}{2}\sum_{1\leq g\leq \lfloor N/2\rfloor} P_2\left(\dfrac{ag}{\gcd(c,N)}\right)r_g.
\end{align*}
\end{remark}

\section{Proofs of identity (\ref{eq14}) and Theorem \ref{thm11}}\label{sec3}

We first derive the level $12$ identity (\ref{eq14}) and prove Theorem \ref{thm11} using sums of generalized $\eta$-quotients on $\Gamma_0(12)$. 

\begin{proof}[Proof of (\ref{eq14})]
We replace $q$ with $q^{12}$ and set $b=q^6$ in (\ref{eq110}) so that

\begin{align}
	\sum_{n=-\infty}^\infty\left[\dfrac{aq^{12n}}{(1-aq^{12n})^2}-\dfrac{q^{12n+6}}{(1-q^{12n+6})^2}\right]=aq^{-3}\Pi_{q^6}^2\dfrac{f^2(-aq^6,-\frac{q^6}{a})}{f^2(-a,-\frac{q^{12}}{a})}. \label{eq31}
\end{align}
Substituting $a=q, q^3, q^5$ in (\ref{eq31}) and adding the resulting identities lead to
\begin{align}
	\dfrac{1}{\Pi_{q^6}^2}\left(\sum_{n\geq 1}\dfrac{q^{2n-1}}{(1-q^{2n-1})^2}-6\sum_{n\geq 1}\dfrac{q^{12n-6}}{(1-q^{12n-6})^2}\right) =: h_1(\tau)+1+\dfrac{1}{h_1(\tau)},\label{eq32}
\end{align}
where 
\begin{align*}
	h_1(\tau) := \dfrac{\eta_{12, 5}^2(\tau)}{\eta_{12, 1}^2(\tau)}
\end{align*}
is modular on $\Gamma_1(12)$ by Lemma \ref{lem27} and satisfies $h_1(\gamma\tau) = 1/h_1(\tau)$ by Lemma \ref{lem26} with $\gamma := [\begin{smallmatrix}
	5 &2\\12 & 5
\end{smallmatrix}]$. Since $\Gamma_1(12)$ and $\gamma$ generate $\Gamma_0(12)$, we see that $j_1(\tau):=h_1(\tau)+1/h_1(\tau)$ is modular on $\Gamma_0(12)$. We next use Remark \ref{rem29} to find the orders of $j_1(\tau)$ at each element of the set
\begin{align}\label{eq33}
	S_{\Gamma_1(12)} = \{\infty, 0,1/5,1/2,1/3,5/3,1/4,7/4,1/6,5/12\}
\end{align}
of inequivalent cusps of $\Gamma_1(12)$ found via Lemma \ref{lem22}. Observe that for such an element $r$, we have
\begin{align*}
	\mathrm{ord}_r\,j_1(\tau) &\geq \min\{\mathrm{ord}_r\, h_1(\tau),\mathrm{ord}_r\, 1/h_1(\tau)\},
\end{align*}
and equality holds when $\mathrm{ord}_r\, h_1(\tau)\neq \mathrm{ord}_r\,1/h_1(\tau)$. We tabulate the values of these orders as shown in Table \ref{tab:tbl31}. 

\begin{table}[h]
	\caption{The orders of $h_1(\tau), 1/h_1(\tau)$, and $j_1(\tau)$ at the cusps of $\Gamma_1(12)$}\label{tab:tbl31}%
	\begin{tabular}{@{}cccc@{}}
		\toprule
		cusp $r$ of $\Gamma_1(12)$ & $\infty$ & $5/12$ & $r\notin\{\infty,5/12\}$ \\
		\midrule
		$\mathrm{ord}_r\, h_1(\tau)$ & $-2$ & $2$ & $0$\\
		\midrule
		$\mathrm{ord}_r\, 1/h_1(\tau)$ & $2$ & $-2$ & $0$\\
		\midrule
		$\mathrm{ord}_r\,j_1(\tau)$ & $-2$ & $-2$ & $\geq 0$ \\
		\bottomrule
	\end{tabular}
\end{table}
Recalling that $j_1(\tau)$ is modular on $\Gamma_0(12)$, we find the set of inequivalent cusps of $\Gamma_0(12)$ given by
\begin{align}\label{eq34}
	S_{\Gamma_0(12)} = \{\infty, 0, 1/2,1/3,1/4,1/6\}
\end{align}
by Lemma \ref{lem21}. We also deduce from that lemma that the following pairs of cusps of $\Gamma_1(12)$ are equivalent under $\Gamma_0(12)$: $(\infty, 5/12), (0,1/5), (1/4,7/4)$, and $(1/3,5/3)$. We see from Table \ref{tab:tbl31} and \cite[Lemma 2]{radu} that the orders of $j_1(\tau)$ at each element of $S_{\Gamma_0(12)}$ are as follows:
\begin{align*}
	\mathrm{ord}_{\infty}\,j_1(\tau) &= \mathrm{ord}_{5/12}\,j_1(\tau) = -2,\\
	\mathrm{ord}_{0}\,j_1(\tau) &= \mathrm{ord}_{1/5}\,j_1(\tau) \geq 0,\\
	\mathrm{ord}_{1/4}\,j_1(\tau) &= \mathrm{ord}_{7/4}\,j_1(\tau) \geq 0,\\
	\mathrm{ord}_{1/3}\,j_1(\tau) &= \mathrm{ord}_{5/3}\,j_1(\tau) \geq 0,\\
	\mathrm{ord}_{1/2}\,j_1(\tau) &\geq 0,\\ 
	\mathrm{ord}_{1/6}\,j_1(\tau) &\geq 0.
\end{align*} 
Thus, $j_1(\tau)\in\mathcal{M}^\infty(12)$ with unique double pole at $\infty$. We note that the $\eta$-quotient
\begin{align*}
	h := \dfrac{\Pi_{q^2}}{\Pi_{q^6}} = \dfrac{\eta^4(4\tau)\eta^2(6\tau)}{\eta^2(2\tau)\eta^4(12\tau)}
\end{align*}
generates $\mathcal{M}^\infty(12)$ with unique simple pole at $\infty$ by Lemmas \ref{lem24} and \ref{lem25}. By Lemma \ref{lem23}, $j_1(\tau)$ is a quadratic polynomial in $h$. From the $q$-expansions 
\begin{align*}
	j_1(\tau) &= q^{-2}+2q^{-1}+3+4q+6q^2+O(q^3),\\
	h &= q^{-1}+2q+q^3+O(q^4),
\end{align*}
we deduce that
\begin{align}\label{eq35}
	j_1(\tau) = h_1(\tau)+\dfrac{1}{h_1(\tau)} = h^2+2h-1.
\end{align}
Comparing (\ref{eq32}) and (\ref{eq35}), and using the definition of $h$, we obtain identity (\ref{eq14}) as desired.
\end{proof}

We require the following result about a modular function on $\mathcal{M}^\infty(12)$ involving 
\begin{align*}
h_2(\tau) := \dfrac{\Pi_{q^6}^2}{\Pi_{q^{12}}^2}=\dfrac{\eta^{12}(12\tau)}{\eta^4(6\tau)\eta^8(24\tau)}.
\end{align*}

\begin{lemma}\label{lem31}
The function 
\begin{align*}
	H:=H(\tau) := h_2(\tau)+\dfrac{16}{h_2(\tau)}
\end{align*}
is modular on $\Gamma_0(12)$, and $hH$ is in $\mathcal{M}^\infty(12)$ with $\mathrm{ord}_\infty\, hH=-4$.
\end{lemma}

\begin{proof}
Observe that by Lemma \ref{lem24}, $h_2(\tau)$ is modular on $\Gamma_0(24)$. Since $\Gamma_0(12)$ is generated by $\Gamma_0(24)$ and $\alpha:=[\begin{smallmatrix}
	1 & 0\\12 & 1
\end{smallmatrix}]$, we know that $H(\tau) = h(\tau)+h(\alpha\tau)$ is modular on $\Gamma_0(12)$. We will show that
\begin{align}\label{eq36}
	h_2(\alpha \tau) = \dfrac{16}{h_2(\tau)}.
\end{align}
Consider the function $H_0(\tau) := h_2(\alpha\tau)h_2(\tau)$, which is modular on $\Gamma_0(24)$ since $\Gamma_0(24)$ is a normal subgroup of $\Gamma_0(12)$. We find the orders of $H_0(\tau)$ at the cusps of $\Gamma_0(24)$ given by 
\begin{align*}
	S_{\Gamma_0(24)} := \{\infty, 0, 1/2,1/3,1/4,1/6,1/8,1/12\}
\end{align*}
from Lemma \ref{lem21} as follows. We first look for the cusps $r_1, r_2\in S_{\Gamma_0(24)}$ such that $\alpha(r_1)$ is equivalent to $r_2$ over $\Gamma_0(24)$, so that by \cite[Lemma 2]{radu}, we have 
\begin{align}\label{eq37}
	\mbox{ord}_{r_1}\,h_2(\alpha\tau)=\mbox{ord}_{\alpha(r_1)}\,h_2(\tau) = \mbox{ord}_{r_2}\,h_2(\tau).
\end{align}
We denote this equivalence by $\alpha(r_1)\sim r_2$. Using Lemma \ref{lem21}, we find that
\begin{align*}
	&\alpha(0) \sim 0, \qquad \alpha(1/2) \sim 1/2, \qquad \alpha(1/3) \sim 1/3, \qquad \alpha(1/4)\sim 1/8,\\
	&\alpha(1/6) \sim 1/6, \qquad \alpha(1/8) \sim 1/4, \qquad \alpha(1/12) \sim \infty,\qquad \alpha(\infty)\sim 1/12.
\end{align*}
We now compute the required orders using Lemma \ref{lem25} and (\ref{eq37}), as shown in Table \ref{tab:tbl32} below. We find from this table that $H_0(\tau)$ is holomorphic at every cusp of $\Gamma_0(24)$, so it is constant by \cite[Lemma 5]{radu}. 

\begin{table}[h]
	\caption{The orders of $h_2(\tau), h_2(\alpha\tau)$ and $H_0(\tau)$ at the cusps of $\Gamma_0(24)$}\label{tab:tbl32}
	\begin{tabular}{@{}cccccc@{}}
		\toprule
		cusp $r$ of $\Gamma_0(24)$ & $\infty$ & $1/12$ & $1/4$ & $1/8$ & $0,1/2,1/3,1/6$\\
		\midrule
		$\mathrm{ord}_r\,h_2(\tau)$ & $-3$ & $3$ & $1$ & $-1$ & $0$\\
		\midrule
		$\mathrm{ord}_r\,h_2(\alpha\tau)$ & $3$ & $-3$ & $-1$ & $1$ & $0$\\
		\midrule
		$\mathrm{ord}_r\,H_0(\tau)$ & $0$ & $0$ & $0$ & $0$ & $0$\\
		\bottomrule
	\end{tabular}
\end{table}

We now employ the transformation formula \cite[Lemma 1]{yang1}
\begin{align*}
	\eta\left(-\dfrac{1}{\tau}\right) = \sqrt{-i\tau}\eta(\tau)
\end{align*}
to compute the constant $H_0(0) = h_2^2(0)$. We obtain 
\begin{align*}
	h_2(0)&=\lim_{\tau\rightarrow\infty}\dfrac{\eta^{12}(-12/\tau)}{\eta^4(-6/\tau)\eta^{8}(-24/\tau)}=4\lim_{\tau\rightarrow\infty}\dfrac{(q^{1/12};q^{1/12})_\infty^{12}}{(q^{1/6};q^{1/6})_\infty^4(q^{1/24};q^{1/24})_\infty^8}=4,
\end{align*}
so that $H_0(0)=16$. Thus, (\ref{eq36}) follows. We now look at Table \ref{tab:tbl32} to get the orders of $H(\tau)$ at the cusps of $\Gamma_0(24)$, as shown in Table \ref{tab:tbl33}.
We use the fact that for such a cusp $r$, we have
\begin{align*}
	\mathrm{ord}_r\,H(\tau) &\geq \min\{\mathrm{ord}_r\, h_2(\tau),\mathrm{ord}_r\, 1/h_2(\tau)\},
\end{align*}
and equality holds when $\mathrm{ord}_r\, h_2(\tau)\neq \mathrm{ord}_r\, 1/h_2(\tau)$.

\begin{table}[h!]
	\caption{The orders of $h_2(\tau), 1/h_2(\tau)$, and $H(\tau)$ at the cusps of $\Gamma_0(24)$}\label{tab:tbl33}
	\begin{tabular}{@{}cccccc@{}}
		\toprule
		cusp $r$ of $\Gamma_0(24)$ & $\infty$ & $1/12$ & $1/4$ & $1/8$ & $0,1/2,1/3,1/6$\\
		\midrule
		$\mathrm{ord}_r\,h_2(\tau)$ & $-3$ & $3$ & $1$ & $-1$ & $0$\\
		\midrule
		$\mathrm{ord}_r\,1/h_2(\tau)$ & $3$ & $-3$ & $-1$ & $1$ & $0$\\
		\midrule
		$\mathrm{ord}_r\,H(\tau)$ & $-3$ & $-3$ & $-1$ & $-1$ & $\geq 0$\\
		\bottomrule
	\end{tabular}
\end{table}

We know that $H(\tau)$ is modular at $\Gamma_0(12)$, whose set $S_{\Gamma_0(12)}$ of its inequivalent cusps is given by (\ref{eq34}). Moreover, we obtain the following pairs of cusps of $\Gamma_0(24)$ that are equivalent under $\Gamma_0(12)$, namely $(\infty, 1/12)$ and $(1/4, 1/8)$. We now deduce the orders of $H(\tau)$ at each element of $S_{\Gamma_0(12)}$ using Table \ref{tab:tbl33} and \cite[Lemma 2]{radu}, with $\infty$ now seen as a cusp of $\Gamma_0(12)$:
\begin{align*}
	\mathrm{ord}_{\infty}\,H(\tau) &= \mathrm{ord}_{1/24}\,H(\tau) = -3,\\
	\mathrm{ord}_{1/4}\,H(\tau) &= \mathrm{ord}_{1/8}\,H(\tau) =-1,\\
	\mathrm{ord}_{0}\,H(\tau) &\geq 0,\\
	\mathrm{ord}_{1/2}\,H(\tau) &\geq 0,\\
	\mathrm{ord}_{1/3}\,H(\tau) &\geq 0,\\
	\mathrm{ord}_{1/6}\,H(\tau) &\geq 0.
\end{align*}
We know from Lemma \ref{lem25} that $\mathrm{ord}_\infty\,h = -1$ and $\mathrm{ord}_{1/4}\,h = 1$. Hence, $hH$ has a unique pole of order $-4$ at $\infty$ and holomorphic elsewhere, so that $hH\in\mathcal{M}^\infty(12)$ as desired.
\end{proof}

\begin{proof}[Proof of Theorem \ref{thm11}]
Setting $a=q^2$ and $a=q^4$ in (\ref{eq31}) and adding the resulting identities yield
\begin{align}\label{eq38}
	\sum_{n\geq 1}\left[\dfrac{q^{2n}}{(1-q^{2n})^2}-\dfrac{q^{6n}}{(1-q^{6n})^2}\right]-4\sum_{n\geq 1}\dfrac{q^{12n-6}}{(1-q^{12n-6})^2}=\Pi_{q^6}^2\left(h_3(\tau)+\dfrac{1}{h_3(\tau)}\right),
\end{align}
where 
\begin{align*}
	h_3(\tau) := \dfrac{\eta_{12, 4}^2(\tau)}{\eta_{12, 2}^2(\tau)}
\end{align*}
is modular on $\Gamma_1(12)$ by Lemma \ref{lem27} and satisfies $h_3(\gamma\tau) = 1/h_3(\tau)$ by Lemma \ref{lem26}. We see that $j_3(\tau) := h_3(\tau)+1/h_3(\tau)$ is modular on $\Gamma_0(12)$. Following Remark \ref{rem29}, we look for the orders of $j_3(\tau)$ at each element of $S_{\Gamma_1(12)}$ given by (\ref{eq33}). Note that for such an element $r$, we have
\begin{align*}
	\mathrm{ord}_r\,j_3(\tau) &\geq \min\{\mathrm{ord}_r\, h_3(\tau),\mathrm{ord}_r\, 1/h_3(\tau)\},
\end{align*}
and equality holds when $\mathrm{ord}_r\, h_3(\tau)\neq \mathrm{ord}_r\,1/h_3(\tau)$. We record these orders as shown in Table \ref{tab:tbl34}. 

\begin{table}[h]
	\caption{The orders of $h_3(\tau), 1/h_3(\tau)$, and $j_3(\tau)$ at the cusps of $\Gamma_1(12)$}\label{tab:tbl34}%
	\begin{tabular}{@{}cccc@{}}
		\toprule
		cusp $r$ of $\Gamma_1(12)$ & $\infty, 5/12$ & $1/4, 7/4$ & $r\notin\{\infty,5/12,1/4,7/4\}$ \\
		\midrule
		$\mathrm{ord}_r\, h_3(\tau)$ & $-1$ & $1$ & $0$\\
		\midrule
		$\mathrm{ord}_r\, 1/h_3(\tau)$ & $1$ & $-1$ & $0$\\
		\midrule
		$\mathrm{ord}_r\,j_3(\tau)$ & $-1$ & $-1$ & $\geq 0$ \\
		\bottomrule
	\end{tabular}
\end{table}

We now compute the orders of $j_3(\tau)$ at each element of $S_{\Gamma_0(12)}$ given by (\ref{eq34}). Recall the following pairs of cusps of $\Gamma_1(12)$ that are equivalent under $\Gamma_0(12)$: $(\infty, 5/12), (0,1/5), (1/4,7/4)$, and $(1/3,5/3)$. From Table \ref{tab:tbl34} and \cite[Lemma 2]{radu}, we find that
\begin{align*}
	\mathrm{ord}_{\infty}\,j_3(\tau) &= \mathrm{ord}_{5/12}\,j_3(\tau) = -1,\\
	\mathrm{ord}_{1/4}\,j_3(\tau) &= \mathrm{ord}_{7/4}\,j_3(\tau) =-1,\\
	\mathrm{ord}_{0}\,j_3(\tau) &= \mathrm{ord}_{1/5}\,j_3(\tau) \geq 0,\\
	\mathrm{ord}_{1/3}\,j_3(\tau) &= \mathrm{ord}_{5/3}\,j_3(\tau) \geq 0,\\
	\mathrm{ord}_{1/2}\,j_3(\tau) &\geq 0,\\ 
	\mathrm{ord}_{1/6}\,j_3(\tau) &\geq 0.
\end{align*} 
Since $\mathrm{ord}_{1/4}\,h=1$, we have that $hj_3(\tau)\in\mathcal{M}^\infty(12)$ with unique double pole at $\infty$, so $hj_3(\tau)$ is a quadratic polynomial in $h$ according to
Lemma \ref{lem23}. By looking at the $q$-expansion
\begin{align*}
	j_3(\tau) = q^{-1}+3q-q^3+O(q^4),
\end{align*}
we obtain 
\begin{align*}
	hj_3(\tau) = h\left(h_3(\tau)+\dfrac{1}{h_3(\tau)}\right) = h^2+1,
\end{align*}
so that the right-hand side of (\ref{eq38}) becomes
\begin{align}\label{eq39}
	\sum_{n\geq 1}\left[\dfrac{q^{2n}}{(1-q^{2n})^2}-\dfrac{q^{6n}}{(1-q^{6n})^2}\right]-4\sum_{n\geq 1}\dfrac{q^{12n-6}}{(1-q^{12n-6})^2}=\Pi_{q^6}^2\left(h+\dfrac{1}{h}\right).
\end{align}
Invoking (\ref{eq14}), the definition of $h$, and the identity
\begin{align}
	\sum_{n\geq 1}\left[\dfrac{q^n}{(1-q^n)^2}-\dfrac{q^{2n}}{(1-q^{2n})^2}\right]=\sum_{n\geq 1}\dfrac{q^{2n-1}}{(1-q^{2n-1})^2},\label{eq310}
\end{align} 
we transform (\ref{eq39}) into
\begin{align}
	\dfrac{1}{\Pi_{q^6}^2}\sum_{n\geq 1}\left[\dfrac{q^n}{(1-q^n)^2}-6\dfrac{q^{6n}}{(1-q^{6n})^2}\right]&-h^2-3h-\dfrac{1}{h}\nonumber\\
	&=\dfrac{5}{\Pi_{q^6}^2}\sum_{n\geq 1}\left[\dfrac{q^{12n-6}}{(1-q^{12n-6})^2}-\dfrac{q^{12n}}{(1-q^{12n})^2}\right].\label{eq311}
\end{align}
We replace $q$ with $q^6$ in (\ref{eq11}), substitute it on the right-hand side of (\ref{eq311}),  and then multiply both sides by $24$. We arrive at 
\begin{align}\label{eq312}
	\dfrac{1}{\Pi_{q^6}^2}\left[24\left(\sum_{n\geq 1}\dfrac{q^n}{(1-q^n)^2}-6\sum_{n\geq 1}\dfrac{q^{6n}}{(1-q^{6n})^2}\right)+5\right]-24h^2-72h-\dfrac{24}{h}=5H.
\end{align}
By Lemma \ref{lem31}, we have $hH\in\mathcal{M}^\infty(12)$ with $\mathrm{ord}_\infty\,hH=-4$. Consequently, by Lemma \ref{lem22}, $hH$ is a quartic polynomial in $h$. It follows from the $q$-expansion
\begin{align*}
	hH = q^{-4}+2q^{-2}+1+20q^2+O(q^4)
\end{align*}
that
\begin{align*}
	hH = h^4-6h^2-3,
\end{align*}
so we conclude from (\ref{eq312}) that
\begin{align*}
	\dfrac{1}{\Pi_{q^6}^2}\left[24\left(\sum_{n\geq 1}\dfrac{q^n}{(1-q^n)^2}-6\sum_{n\geq 1}\dfrac{q^{6n}}{(1-q^{6n})^2}\right)+5\right]=5h^3+24h^2+42h+\dfrac{9}{h},
\end{align*}
which is precisely (\ref{eq17}). This completes the proof.
\end{proof}

\section{Proof of Theorem \ref{thm12}}\label{sec4}

We next prove Theorem \ref{thm12} by employing certain sums of generalized $\eta$-quotients on $\Gamma_0(16)$.

\begin{proof}[Proof of (\ref{eq18})]
Replacing $q$ with $q^{16}$ and setting $b=q^8$ in (\ref{eq110}) give
\begin{align}
	\sum_{n=-\infty}^\infty\left[\dfrac{aq^{16n}}{(1-aq^{16n})^2}-\dfrac{q^{16n+8}}{(1-q^{16n+8})^2}\right]=aq^{-4}\Pi_{q^8}^2\dfrac{f^2(-aq^8,-\frac{q^8}{a})}{f^2(-a,-\frac{q^{16}}{a})}. \label{eq41}
\end{align}	
Substituting $a=q, q^3, q^5, q^7$ in (\ref{eq41}) and combining the resulting identities, we have 
{\small \begin{align}
		\dfrac{1}{\Pi_{q^8}^2}\left(\sum_{n\geq 1}\dfrac{q^{2n-1}}{(1-q^{2n-1})^2}-8\sum_{n\geq 1}\dfrac{q^{16n-8}}{(1-q^{16n-8})^2}\right) =: s_1(\tau)+s_2(\tau)+\dfrac{1}{s_1(\tau)}+\dfrac{1}{s_2(\tau)},\label{eq42}
\end{align}}%
where 
\begin{align*}
	s_1(\tau) := \dfrac{\eta_{16, 7}^2(\tau)}{\eta_{16, 1}^2(\tau)}, \qquad s_2(\tau) := \dfrac{\eta_{16, 5}^2(\tau)}{\eta_{16, 3}^2(\tau)}
\end{align*}
are modular on $\Gamma_1(16)$ by Lemma \ref{lem27}. By Lemma \ref{lem26}, we have
\begin{align*}
	s_1(\gamma'\tau) = s_2(\tau), \qquad s_2(\gamma'\tau) = \dfrac{1}{s_1(\tau)},
\end{align*}
where $\gamma' := [\begin{smallmatrix}
	3 &-1\\16 & -5
\end{smallmatrix}]$. As $\Gamma_1(16)$ and $\gamma'$ generate $\Gamma_0(16)$, we deduce that 
\begin{align*}
	z(\tau) := s_1(\tau)+s_2(\tau)+\dfrac{1}{s_1(\tau)}+\dfrac{1}{s_2(\tau)}
\end{align*}
is modular on $\Gamma_0(16)$. In view of Remark \ref{rem29}, we compute the orders of $z(\tau)$ at each element of the set
\begin{align}\label{eq43}
	S_{\Gamma_1(16)} = \{\infty,0,1/3,1/5,1/7,1/2,1/6,1/4,3/4,1/8,3/8,3/16,5/16,7/16\}
\end{align}
of inequivalent cusps of $\Gamma_1(16)$ obtained via Lemma \ref{lem22}. For such an element $r$, we know that
\begin{align*}
	\mathrm{ord}_r\,z(\tau) &\geq \min\{\mathrm{ord}_r\, s_1(\tau),\mathrm{ord}_r\, s_2(\tau),\mathrm{ord}_r\, 1/s_1(\tau),\mathrm{ord}_r\, s_2(\tau)\},
\end{align*}
and equality holds when $\mathrm{ord}_r\, s_1(\tau),\mathrm{ord}_r\, s_2(\tau),\mathrm{ord}_r\, 1/s_1(\tau)$, and $\mathrm{ord}_r\, s_2(\tau)$ are not all equal. We record these orders as shown in Table \ref{tab:tbl41}. 

\begin{table}[h]
	\caption{The orders of $s_j(\tau), 1/s_j(\tau)$ (for $j\in\{1,2\}$), and $z(\tau)$ at the cusps of $\Gamma_1(16)$}\label{tab:tbl41}%
	\begin{tabular}{@{}cccccc@{}}
		\toprule
		cusp $r$ of $\Gamma_1(16)$ & $\infty$ & $3/16$ & $5/16$ & $7/16$ &$r\notin\{\infty,3/16,5/16,7/16\}$\\
		\midrule
		$\mathrm{ord}_r\, s_1(\tau)$ & $-3$ & $-1$ & $1$ & $3$ & $0$\\
		\midrule
		$\mathrm{ord}_r\, s_2(\tau)$ & $-1$ & $3$ & $-3$ & $1$ & $0$\\
		\midrule
		$\mathrm{ord}_r\, 1/s_1(\tau)$ & $3$ & $1$ & $-1$ & $-3$ & $0$\\
		\midrule
		$\mathrm{ord}_r\, 1/s_2(\tau)$ & $1$ & $-3$ & $3$ & $-1$ & $0$\\
		\midrule
		$\mathrm{ord}_r\,z(\tau)$ & $-3$ & $-3$ & $-3$ & $-3$ & $\geq 0$ \\
		\bottomrule
	\end{tabular}
\end{table}

Recall that $z(\tau)$ is modular on $\Gamma_0(16)$. By Lemma \ref{lem21}, we have the set of inequivalent cusps of $\Gamma_0(16)$ given by
\begin{align}\label{eq44}
	S_{\Gamma_0(16)} = \{\infty, 0, 1/2,1/4,3/4,1/8\}.
\end{align}
We also enumerate the following pairs and quadruples of the cusps of $\Gamma_1(16)$ that are equivalent under $\Gamma_0(16)$: $(1/2,1/6), (1/8,3/8), (0,1/3,1/5,1/7)$ and $(\infty,3/16,5/16,7/16)$. We infer from Table \ref{tab:tbl41} and \cite[Lemma 2]{radu} that the orders of $z(\tau)$ at each element of $S_{\Gamma_0(16)}$ are as follows:
\begin{align*}
	\mathrm{ord}_{\infty}\,z(\tau) &= \mathrm{ord}_{3/16}\,z(\tau) = \mathrm{ord}_{5/16}\,z(\tau)=  \mathrm{ord}_{7/16}\,z(\tau)=-3,\\
	\mathrm{ord}_{0}\,z(\tau) &= \mathrm{ord}_{1/3}\,z(\tau) = \mathrm{ord}_{1/5}\,z(\tau)=  \mathrm{ord}_{1/7}\,z(\tau)\geq 0,\\
	\mathrm{ord}_{1/2}\,z(\tau) &= \mathrm{ord}_{1/6}\,z(\tau) \geq 0,\\
	\mathrm{ord}_{1/8}\,z(\tau) &= \mathrm{ord}_{3/8}\,z(\tau) \geq 0,\\
	\mathrm{ord}_{1/4}\,z(\tau) &\geq 0,\\ 
	\mathrm{ord}_{3/4}\,z(\tau) &\geq 0.
\end{align*}  
We see that $z(\tau)\in\mathcal{M}^\infty(16)$ with unique triple pole at $\infty$. As the $\eta$-quotient 
\begin{align*}
	s := \dfrac{\Pi_{q^4}}{\Pi_{q^8}}=\dfrac{\eta^6(8\tau)}{\eta^2(4\tau)\eta^4(16\tau)}
\end{align*}
generates $\mathcal{M}^\infty(16)$ with unique simple pole at $\infty$ by Lemmas \ref{lem24} and \ref{lem25}, we know from Lemma \ref{lem23} that $z(\tau)$ is a cubic polynomial in $s$. We imply from 
the $q$-expansion 
\begin{align*}
	z(\tau) = q^{-3}+2q^{-2}+4q^{-1}+4+6q+8q^2+O(q^3)
\end{align*}
that 
\begin{align}\label{eq45}
	z(\tau) = s^3+2s^2+4s+4.
\end{align}
Comparing (\ref{eq42}) and (\ref{eq45}), and using the definition of $s$, lead us to (\ref{eq18}) as desired. 
\end{proof}

To establish (\ref{eq19}), we require the following lemma concerning a modular function on $\Gamma_0(16)$ that involves
\begin{align*}
s_3(\tau) := \dfrac{\Pi_{q^8}^2}{\Pi_{q^{16}}^2} = \dfrac{\eta^{12}(16\tau)}{\eta^4(8\tau)\eta^8(32\tau)}.
\end{align*}

\begin{lemma}\label{lem41}
The function 
\begin{align*}
	S := S(\tau) = s_3(\tau)+\dfrac{16}{s_3(\tau)}
\end{align*}
lies on $\mathcal{M}^\infty(16)$ with $\mathrm{ord}_\infty\,S = -4$.
\end{lemma}

\begin{proof}
We know from Lemma \ref{lem24} that $s_3(\tau)$ is modular on $\Gamma_0(32)$. Since $\Gamma_0(32)$ and $\beta:= [\begin{smallmatrix}
	1 & 0\\16 &1 
\end{smallmatrix}]$ generate $\Gamma_0(16)$, we deduce that $S(\tau) = s_3(\tau)+s_3(\beta\tau)$ is modular on $\Gamma_0(16)$. Following similar arguments as in the proof of (\ref{eq36}) from Lemma \ref{lem31}, we get 
\begin{align*}
	s_3(\beta\tau) = \dfrac{16}{s_3(\tau)}.
\end{align*}
We next apply Lemma \ref{lem21} and obtain the set 
\begin{align*}
	S_{\Gamma_0(32)} = \{\infty, 1/2, 1/4, 3/4, 1/8, 3/8, 1/16\}
\end{align*}
of inequivalent cusps of $\Gamma_0(32)$.  We now use Lemma \ref{lem25} to compute the orders of $S(\tau)$ at the cusps of $\Gamma_0(32)$, as shown in Table \ref{tab:tbl42}.
We remark that for such a cusp $r$, we have
\begin{align*}
	\mathrm{ord}_r\,S(\tau) &\geq \min\{\mathrm{ord}_r\, s_3(\tau),\mathrm{ord}_r\, 1/s_3(\tau)\},
\end{align*}
and equality holds when $\mathrm{ord}_r\, s_3(\tau)\neq \mathrm{ord}_r\, 1/s_3(\tau)$.

\begin{table}[h]
	\caption{The orders of $s_3(\tau), 1/s_3(\tau)$, and $S(\tau)$ at the cusps of $\Gamma_0(32)$}\label{tab:tbl42}
	\begin{tabular}{@{}cccc@{}}
		\toprule
		cusp $r$ of $\Gamma_0(32)$ & $\infty$ & $1/16$ & $r\notin\{\infty,1/16\}$\\
		\midrule
		$\mathrm{ord}_r\,s_3(\tau)$ & $-4$ & $4$ & $0$ \\
		\midrule
		$\mathrm{ord}_r\,1/s_3(\tau)$ & $4$ & $-4$ & $0$ \\
		\midrule
		$\mathrm{ord}_r\,S(\tau)$ & $-4$ & $-4$ & $\geq 0$ \\
		\bottomrule
	\end{tabular}
\end{table}

Recall that $S(\tau)$ is modular on $\Gamma_0(16)$, whose set $S_{\Gamma_0(16)}$ of its inequivalent cusps is given by (\ref{eq44}). In addition, we find $(\infty, 1/16)$ and $(1/8, 3/8)$ as the only pairs of cusps of $\Gamma_0(32)$ that are equivalent under $\Gamma_0(16)$. We now deduce the orders of $S(\tau)$ at each element of $S_{\Gamma_0(16)}$ using Table \ref{tab:tbl42} and \cite[Lemma 2]{radu}, with $\infty$ now considered as a cusp of $\Gamma_0(16)$:
\begin{align*}
	\mathrm{ord}_{\infty}\,S(\tau) &= \mathrm{ord}_{1/32}\,S(\tau) = -4,\\
	\mathrm{ord}_{1/8}\,S(\tau) &= \mathrm{ord}_{3/8}\,S(\tau) \geq 0,\\
	\mathrm{ord}_{0}\,S(\tau) &\geq 0,\\
	\mathrm{ord}_{1/2}\,S(\tau) &\geq 0,\\
	\mathrm{ord}_{1/4}\,S(\tau) &\geq 0,\\
	\mathrm{ord}_{3/4}\,S(\tau) &\geq 0.
\end{align*}
Hence, we conclude that $S\in \mathcal{M}^\infty(16)$ with $\mathrm{ord}_\infty\, S=-4$ as desired.
\end{proof}

\begin{proof}[Proof of (\ref{eq19})]
Setting $a=q^2, q^4, q^6$ in (\ref{eq41}) and combining the resulting identities yield
\begin{align}\label{eq46}
	\sum_{n\geq 1}\left[\dfrac{q^{2n}}{(1-q^{2n})^2}-\dfrac{q^{8n}}{(1-q^{8n})^2}\right]-6\sum_{n\geq 1}\dfrac{q^{16n-8}}{(1-q^{16n-8})^2}=\Pi_{q^8}^2\left(s_4(\tau)+1+\dfrac{1}{s_4(\tau)}\right),
\end{align}
where 
\begin{align*}
	s_4(\tau) := \dfrac{\eta_{16, 6}^2(\tau)}{\eta_{16, 2}^2(\tau)}
\end{align*}
is modular on $\Gamma_1(16)$ by Lemma \ref{lem27} and satisfies $s_4(\gamma'\tau) = 1/s_4(\tau)$ by Lemma \ref{lem26}. Thus, $j_4(\tau) := s_4(\tau)+1/s_4(\tau)$ is modular on $\Gamma_0(16)$. Appealing to Remark \ref{rem29}, we seek for the orders of $j_4(\tau)$ at each element of $S_{\Gamma_1(16)}$ given by (\ref{eq43}). For such an element $r$, we have that
\begin{align*}
	\mathrm{ord}_r\,j_4(\tau) &\geq \min\{\mathrm{ord}_r\, s_4(\tau),\mathrm{ord}_r\, 1/s_4(\tau)\},
\end{align*}
and equality holds when $\mathrm{ord}_r\, s_4(\tau)\neq \mathrm{ord}_r\,1/s_4(\tau)$. We present these orders as shown in Table \ref{tab:tbl43}. 

\begin{table}[h]
	\caption{The orders of $s_4(\tau), 1/s_4(\tau)$, and $j_4(\tau)$ at the cusps of $\Gamma_1(16)$}\label{tab:tbl43}%
	\begin{tabular}{@{}cccccc@{}}
		\toprule
		cusp $r$ of $\Gamma_1(16)$ & $\infty$ & $3/16$ & $5/16$ & $7/16$ &$r\notin\{\infty,3/16,5/16,7/16\}$\\
		\midrule
		$\mathrm{ord}_r\, s_4(\tau)$ & $-2$ & $2$ & $2$ & $-2$ & $0$\\
		\midrule
		$\mathrm{ord}_r\, 1/s_4(\tau)$ & $2$ & $-2$ & $-2$ & $2$ & $0$\\
		\midrule
		$\mathrm{ord}_r\,j_4(\tau)$ & $-2$ & $-2$ & $-2$ & $-2$ & $\geq 0$ \\
		\bottomrule
	\end{tabular}
\end{table}

We now calculate the orders of $j_4(\tau)$ at each element of $S_{\Gamma_0(16)}$ given by (\ref{eq44}). Recall the following pairs and quadruples of the cusps of $\Gamma_1(16)$ that are equivalent under $\Gamma_0(16)$: $(1/2,1/6), (1/8,3/8)$,\\ $(0,1/3,1/5,1/7)$ and $(\infty,3/16,5/16,7/16)$. From Table \ref{tab:tbl43} and \cite[Lemma 2]{radu}, we see that
\begin{align*}
	\mathrm{ord}_{\infty}\,j_4(\tau) &= \mathrm{ord}_{3/16}\,j_4(\tau) = \mathrm{ord}_{5/16}\,j_4(\tau)=  \mathrm{ord}_{7/16}\,j_4(\tau)=-2,\\
	\mathrm{ord}_{0}\,j_4(\tau) &= \mathrm{ord}_{1/3}\,j_4(\tau) = \mathrm{ord}_{1/5}\,j_4(\tau)=  \mathrm{ord}_{1/7}\,j_4(\tau)\geq 0,\\
	\mathrm{ord}_{1/2}\,j_4(\tau) &= \mathrm{ord}_{1/6}\,j_4(\tau) \geq 0,\\
	\mathrm{ord}_{1/8}\,j_4(\tau) &= \mathrm{ord}_{3/8}\,j_4(\tau) \geq 0,\\
	\mathrm{ord}_{1/4}\,j_4(\tau) &\geq 0,\\ 
	\mathrm{ord}_{3/4}\,j_4(\tau) &\geq 0.
\end{align*}

Thus, $j_4(\tau)\in\mathcal{M}^\infty(16)$ with unique double pole at $\infty$, so it is a quadratic polynomial in $s$ by Lemma \ref{lem23}. We infer from the $q$-expansion
\begin{align*}
	j_4(\tau) = q^{-2}+2+4q^2+O(q^6)
\end{align*}
that
\begin{align*}
	j_4(\tau) = s_4(\tau)+\dfrac{1}{s_4(\tau)} = s^2+2,
\end{align*}  
so the right-hand side of (\ref{eq46}) becomes
\begin{align}\label{eq47}
	\sum_{n\geq 1}\left[\dfrac{q^{2n}}{(1-q^{2n})^2}-\dfrac{q^{8n}}{(1-q^{8n})^2}\right]-6\sum_{n\geq 1}\dfrac{q^{16n-8}}{(1-q^{16n-8})^2}=\Pi_{q^8}^2\left(s^2+3\right).
\end{align}
Invoking (\ref{eq18}), (\ref{eq310}), and the definition of $s$, we write (\ref{eq47}) as
\begin{align}
	\dfrac{1}{\Pi_{q^8}^2}\sum_{n\geq 1}\left[\dfrac{q^n}{(1-q^n)^2}-8\dfrac{q^{8n}}{(1-q^{8n})^2}\right]&-s^3-3s^2-4s-7\nonumber\\
	&=\dfrac{7}{\Pi_{q^8}^2}\sum_{n\geq 1}\left[\dfrac{q^{16n-8}}{(1-q^{16n-8})^2}-\dfrac{q^{16n}}{(1-q^{16n})^2}\right].\label{eq48}
\end{align}
We now replace $q$ with $q^8$ in (\ref{eq11}), put it on the right-hand side of (\ref{eq48}), and then multiply both sides by $24$. We find that
\begin{align}\label{eq49}
	\dfrac{1}{\Pi_{q^8}^2}\left[24\left(\sum_{n\geq 1}\dfrac{q^n}{(1-q^n)^2}-8\sum_{n\geq 1}\dfrac{q^{8n}}{(1-q^{8n})^2}\right)+7\right]-24s^3-72s^2-96s-168=7S.
\end{align}
We know from Lemma \ref{lem41} that $S\in \mathcal{M}^\infty(16)$ with unique pole of order $-4$ at $\infty$. Thus, $S$ is a quartic polynomial in $s$ by Lemma \ref{lem23}. In view of
the $q$-expansion 
\begin{align*}
	S = q^{-4}+20q^4-62q^{12}+O(q^{13}),
\end{align*}
we get 
\begin{align*}
	S = s^4-8,
\end{align*}
and (\ref{eq49}) becomes 
\begin{align*}
	\dfrac{1}{\Pi_{q^8}^2}\left[24\left(\sum_{n\geq 1}\dfrac{q^n}{(1-q^n)^2}-8\sum_{n\geq 1}\dfrac{q^{8n}}{(1-q^{8n})^2}\right)+7\right]=7s^4+24s^3+72s^2+96s+112,
\end{align*}
which is exactly (\ref{eq19}).
\end{proof}

\section*{Acknowledgments}
The author would like to express his gratitude to Dr. M. V. Yathirajsharma for sending him a copy of the paper \cite{yat}. 

\bibliographystyle{amsplain}
\bibliography{gosper1216}
\end{document}